# The robust single machine scheduling problem with uncertain release and processing times

Nitish Umang *  Alan L. Erera †  Michel Bierlaire ‡

*CIRRELT and DIRO, University of Montreal, 2900 Boulevard Edouard-Montpetit, Montreal, QC H3T 1J4, Canada, nitish.umang@cirrelt.ca

†School of Industrial and Systems Engineering, Georgia Institute of Technology, Atlanta, Georgia 30332, USA, alerera@isye.gatech.edu

‡Transport and Mobility Laboratory (TRANSP-OR), School of Architecture, Civil and Environmental Engineering (ENAC), Ecole Polytechnique Fédérale de Lausanne (EPFL), CH-1015 Lausanne, Switzerland, michel.bierlaire@epfl.ch

**Abstract**

In this work, we study the single machine scheduling problem with uncertain release times and processing times of jobs. We adopt a robust scheduling approach, in which the measure of robustness to be minimized for a given sequence of jobs is the worst-case objective function value from the set of all possible realizations of release and processing times. The objective function value is the total flow time of all jobs. We discuss some important properties of robust schedules for zero and non-zero release times, and illustrate the added complexity in robust scheduling given non-zero release times. We propose heuristics based on variable neighborhood search and iterated local search to solve the problem and generate robust schedules. The algorithms are tested and their solution performance is compared with optimal solutions or lower bounds through numerical experiments based on synthetic data.

# 1 Introduction

Scheduling involves the optimal allocation of scarce resources to activities over time. Scheduling problems are an integral part of planning in areas such as production, service, manufacturing and transportation. In the past few decades, the practical importance and complexity of the general scheduling problem has motivated a significant volume of research in a wide variety of scheduling environments including production, manufacturing, transportation and logistics systems. Using standard notation, scheduling problems include a set of $n$ jobs that must be scheduled on a set of $m$ machines subject to certain constraints to optimize a desired objective function. In reality one or more characteristics of the jobs may be uncertain due to factors such as worker performance variability, changes in the work environment, variability in tool quality, and a variety of other factors. In this paper we study the most common configuration i.e. single machine scheduling using $m = 1$, with a particular focus on generating "robust" schedules. Our primary goal is to demonstrate the challenge of building robustness into scheduling solutions, while keeping the problems simple enough to permit useful analysis.

The major emphasis in past scheduling research has been on deterministic problems in which the schedule is computed and fixed in advance assuming perfect knowledge of job-specific attributes such as release times, processing times and/or due dates. However, a major drawback of precomputed schedules is that even small deviations in job parameter values can disrupt the schedule and lead to significant system performance degradation. Thus it is desirable to generate schedules that are "robust" given task parameter uncertainty. Consider a schedule that is created off-line and then placed into operation. During its execution, a disturbance may render the planned schedule infeasible. In response to the disturbance, a control action is executed to restore feasibility. A robust schedule is an *a priori* schedule which maintains high system performance in the presence of stochastic disturbances given a policy for control actions. In this study, we use a simple control policy that shifts the disrupted schedule in time without altering the original planned sequence of jobs, which is particularly useful in situations where changing the sequence may result in additional cost.

In classical stochastic scheduling, uncertain job attributes are modeled as independent random variables with known distributions. The performance of a schedule is dependent on the specific realization of each uncertain parameter during execution, while the design objective typically is to optimize the expected performance of the system. There are drawbacks of this approach. First, it assumes knowledge of probability distributions for the uncertain parameters, which are often unknown and almost never precisely known and may be difficult to estimate. Moreover, the decision maker may be more interested in hedging against the worst-case performance of the system than optimizing the system performance averaged over all possible realizations. However classical approaches fail to recognize this fact.

In this work, we study robust scheduling to determine a schedule which has the best worst-case performance. Our focus is single machine scheduling where the performance criterion is the total flow time of all jobs. The rest of the paper is organized as follows: Section 2 provides a brief literature review on the general scheduling problem with a particular focus on research work done in dealing with uncertainty in the context of machine scheduling problems. In Section 3, we formally define the framework of the robust single machine scheduling problem and provide some important insights into the deterministic and stochastic variants of the single machine scheduling problem. In Section 4 we propose solution algorithms to obtain good solutions for the robust single machine scheduling problem with release times. In Section 5, we present computational results based on artificial instances to test and validate the efficiency of the proposed algorithms and compare their solution performance from a computational perspective. Finally we give some concluding remarks in Section 6.

# 2 Literature Review

Comprehensive literature surveys on the general scheduling problem in a wide variety of scheduling environments can be found in Lawler (1976), Graham et al. (1979) and Blazewicz (1987). Graham et al. (1979) established a three-field notation $\alpha|\beta|\gamma$ to simplify the categorization of different types of machine scheduling problems. In this notation, the parameters $\alpha$, $\beta$ and $\gamma$ describe the machine environment, the job characteristics and the optimality criterion respectively. For example, $1|r_j|\sum C_j$ denotes the variant of the problem in which there is a single machine, each job $j$ is available for processing only at the release time $r_j$ or later, and the objective is to minimize the sum of completion times of all jobs as given by $\sum C_j$. As another example, $1|r_j, \texttt{prec}|\sum_j C_j\text{-}r_j$ denotes the problem of scheduling the jobs with precedence constraints and release times on a single machine with the objective to minimize the total flow times of all jobs. As it will be impossible to enumerate all the variants of the problem and out of the scope of this study, we refer to Graham et al. (1979) for a survey on the different types of scheduling problems in literature.

Research has addressed machine scheduling problems in which one or more aspects of the jobs such as release times, processing times and other job-related properties are random, or the machines are subject to random breakdowns, or both. Glazebrook (1979), Weiss and Pinedo (1980), Emmons and Pinedo (1990) are few examples of such works. Stochastic machine scheduling problems focusing on probabilistic times have been studied by Wu and Zhou (2008), Skutella and Uetz (2005), Cai and Zhou (2005) and Soroush and Fredendall (1994) in which the job attributes are modeled as independent random variables with given distributions, whose actual values are realized during the execution of the schedule after a scheduling decision has been made. Dynamic scheduling methods in which jobs are dispatched dynamically to account for random disruptions in real time are studied by Gittins and Glazebrook (1977), Pinedo (1983), Glazebrook (1981), Glazebrook (1985) and few others. Another line of research focuses on responding to random disruptions that occur in real time, making it impossible to adhere to the originally planned schedule. Bean et al. (1987) and Roundy et al. (1989) are examples of such works. For detailed literature surveys on fundamental approaches for scheduling under uncertainty, refer to Herroelen and Leus (2005), Mohring et al. (1985), Mohring et al. (1984) and Pinedo and Schrage (1982).

Kouvelis and Yu (1997) developed robust versions of many traditional discrete optimization problems. In general three different measures of robustness can be defined; one that minimizes the maximum absolute cost over the set of possible outcomes, a second that minimizes the maximum regret i.e. the absolute difference in the solution cost between the realized outcome and the corresponding optimal solution for the outcome, and a third that minimizes the maximum relative deviation of the realized outcome from the corresponding optimal solution. Daniels and Kouvelis (1995) study the robust single machine scheduling problem without release times in which schedule robustness is measured by the absolute or relative deviation of the realized cost from optimality. They describe properties of robust schedules which allow the selection of a finite set of scenarios from uncertainty intervals of processing times to determine the worst-case deviation from optimality for a given schedule, and propose exact and heuristic solution approaches to obtain robust schedules. Yang and Yu (2002) study the same problem as Daniels and Kouvelis (1995), show that the problem is NP-hard even in the case of two scenarios for all three measures of robustness described earlier, and propose two alternative heuristic methods to obtain robust schedules. Kasperski (2005) studies the single machine scheduling problem for the absolute deviation measure of robustness, the maximum lateness performance criterion, and uncertainty intervals for the processing times. A polynomial time algorithm is proposed to solve the problem. Some other references that study the robust version of the single machine scheduling problem considering uncertainty intervals for processing times are Lebedev and Averbakh (2006), Montemanni (2007) and Kasperski and Zielinski (2008). More recently, Lu et al. (2012) study the single machine scheduling problem with uncertainty in the job processing times and sequence-dependent family setup times. In their study, the performance criterion is the total flow time of jobs, and the measure of schedule robustness is the maximum absolute deviation from the optimal solution in the worst-case scenario. They reformulate the problem as a robust constrained shortest path problem and propose a simulated annealing-based algorithm to determine robust schedules.

In this research, we use the maximum absolute cost over the set of all possible outcomes as the measure of robustness and the total flow time of jobs as the performance criterion to create robust schedules in the context of the single machine scheduling problem. To the best of our knowledge, this is the first paper that considers uncertainty in both release times and processing times in the robust scheduling context for the single machine scheduling problem. We discuss some important properties of robust schedules with zero and non-zero release times, demonstrate the added complexity when non-zero release times are considered, propose an exact method to instantaneously solve the deterministic variant of the single machine scheduling problem with release times, and develop heuristic methods based on variable neighborhood search and iterated local search to generate robust schedules. The solution performance of the proposed algorithms are tested and validated through extensive numerical experiments based on artificial data.

# 3 Robust Single Machine Scheduling Problem

## 3.1 Problem Definition

We consider a set of $n$ jobs that are required to be scheduled on a single machine. We are interested in generating robust schedules for uncertain scheduling environments, in which there is stochastic variability in the release times $r_i$ and the processing times $p_i$ of jobs. In our problem, the release times and the processing times of the jobs are specified as independent ranges of values with unknown probability distributions, such that the release time interval of job $i$ is $[\underline{r_i}, \overline{r_i}]$ and the processing time interval of job $i$ is $[\underline{p_i}, \overline{p_i}]$. Let the infinite set of possible realizations of release times and processing times be represented by the set $\Omega$. Then a possible outcome $\omega \in \Omega$, represents a unique set of release times and processing times of the jobs, that can be realized with a certain positive and unknown probability. Let the decision space consisting of all possible job sequences be given by the set $P$. The cost of making sequencing decision $\sigma \in P$ under scenario $\omega \in \Omega$ is given by $f(\sigma, \omega)$. The optimal decision and the optimal cost under scenario $\omega \in \Omega$ are given by $\sigma_{\omega^*}$ and $f^*(\omega)$ respectively.

We assume the following input data to be available for the singe machine scheduling problem :

| | |
|---|---|
| $N =$ | set of jobs |
| $i =$ | $1, ..., \|N\|$ jobs |
| $\Omega =$ | the infinite set of possible realizations |
| $P =$ | decision space representing the set of all possible sequences |
| $r_i^\omega =$ | release time of job $i \in N$ for the realization $\omega \in \Omega$ |
| $p_i^\omega =$ | processing time of job $i \in N$ for the realization $\omega \in \Omega$ |

The objective in the absolute robust single machine scheduling problem (ARSMSP), can be mathematically expressed as follows

$$(ARSMSP) \min_{\sigma \in P} \{\max_{\omega \in \Omega} (f(\sigma, \omega)\} \tag{1}$$

Let $N_\sigma$ represent the ordered sequence of jobs for the sequence $\sigma \in P$, such that for jobs $i, j \in N_\sigma$ and $j > i$, it is implied that job $j$ is sequenced after job $i$ in $\sigma$. For a given sequence $\sigma \in P$, realization $\omega \in \Omega$ and the performance criterion as the total flow time of jobs, we have

$$f(\sigma, \omega) = \sum_{i \in N_\sigma} (s_i - r_i^\omega + p_i^\omega) \tag{2}$$

subject to the conditions

$$s_1 = r_1^\omega \tag{3}$$

$$s_i = \max\left(r_i^\omega, s_{i-1} + p_{i-1}^\omega\right) \quad \forall i \in N_\sigma, i \geq 2 \tag{4}$$

The deterministic single machine scheduling problem (DSMSP) to determine $f^*(\omega)$ for a given realization $\omega \in \Omega$ can be formulated as follows:

$$(DSMSP) \quad \min \textstyle\sum_{i \in N} (s_i - r_i^\omega + p_i^\omega) \tag{5}$$

$$\text{s.t. } s_i - r_i^\omega \geq 0 \qquad \forall i \in N \tag{6}$$

$$s_j \geq s_i + p_i^\omega || s_i \geq s_j + p_j^\omega \quad \forall i, j \in N, i \neq j \tag{7}$$

In the above formulation, constraints (6) ensure that the processing of a job starts only at or after the release time of the job. Constraints (7) are the disjunctive constraints that ensure that two jobs are not processed at the same time. Unfortunately the disjunctive constraints are non-linear, but can be linearized using the bigM approach, and reformulated as

$$s_j + M(1 - z_{ij}) \geq s_i + p_i^{\omega} \quad \forall i, j \epsilon N, i \neq j \tag{8}$$

$$z_{ij} + z_{ji} = 1 \quad \forall i, j \epsilon N, i \neq j \tag{9}$$

where $z_{ij}$ is a binary variable equal to 1 if job $i$ preceeds job $j$ without overlapping, 0 otherwise, and $M$ is a large positive constant. With regard to complexity, DSMSP is strongly NP-hard (Lenstra et al. (1977)).

In the following section, our aim is to discuss some of the most important results related to the deterministic and robust variants of the single machine scheduling problem, and demonstrate the added complexity when there is uncertainty in both the release times and the processing times of the jobs. We begin by briefly looking at the deterministic version of the single machine scheduling problem without release times.

## 3.2 Scheduling without release times

### 3.2.1 Deterministic Problem

The simplest scheduling problem arises when the release times of all jobs are equal to zero. The obvious approach to solve this problem is to assign a priority to each job based on the optimality criterion, and assign the jobs in the order of decreasing priorities whenever the machine becomes available. Note that in the absence of release times, the flow time of a given job is equivalent to it's completion time. Thus according to the notation discussed earlier, the single machine scheduling problem without release times with the objective to minimize the total flow times can be represented by $1|C_j$. Intuitively, it makes sense to schedule the job with the shortest processing time at the beginning so that the delays to all the other jobs are minimized, and in a similar way, schedule the remaining jobs in the order of increasing processing times. In the literature, this is commonly known as the *Shortest Processing Time* (SPT) rule. We have the following useful result(Smith (1956)).

RESULT 1: *SPT rule is an exact algorithm to solve* $1|\sum C_j$ *with time complexity* $O(n \log n)$*.*

### 3.2.2 Properties of robust schedules without release times

In the following discussion, we discuss some properties of robust schedules with the performance criterion as the total flow time or completion time (both are equivalent for zero release times) of the jobs. The release time of each job $i \in N$ is equal to zero, and the processing time interval of job $i$ is $[\underline{p_i}, \overline{p_i}]$.

**ARSMSP without release times** We begin with a simple result for the absolute robust single machine scheduling problem (ARSMSP) without release times.

RESULT 2: *The optimal solution to the ARSMSP without release times is the sequence of jobs obtained by arranging the jobs in increasing order of* $\overline{p_i}$*, that is the highest processing time values for all jobs.*

Proof: Let the sequence of jobs obtained by arranging the jobs in increasing order of the highest processing times be $\sigma_{\omega_{max}}$. The worst case contingency for this sequence corresponds to the case when each job $i \in N$ assumes its highest processing time $\overline{p_i}$. However it is obvious that the sequence $\sigma_{\omega_{max}}$ is also the optimal decision for the realization corresponding to this worst case contingency (using SPT algorithm discussed earlier). Hence for any other sequence $\sigma \in P$, the flow time for the worst case contingency corresponding to $p = \overline{p_i}$ for each job $i \in N$, is higher than for the sequence $\sigma_{\omega_{max}}$.

## 3.3 Scheduling with release times

### 3.3.1 Deterministic Problem

As stated earlier, the deterministic single machine scheduling problem (DSMSP) is an NP-complete problem. The MILP reformulation of DSMSP is inefficient and slow, and cannot be used to solve large problem size to optimality in a reasonable computational time. However in order to obtain a lower bound to the absolute robust single machine scheduling problem (ARSMSP), it is desirable to have an efficient algorithm to obtain the optimal solution or at the very least a tight upper bound to the DSMSP. This point is better illustrated by the following result.

RESULT 3: *The maximum optimal value* $f^*(\omega)$ *over the set of all possible realizations* $\omega \in \Omega$ *is a lower bound to the absolute robust single machine scheduling problem (ARSMSP) with (or without) release times.*

Proof: Let's say that we are given a sequence $\sigma \in P$, for which $\omega_\sigma$ is the worst case realization. Then we have

$$f(\sigma, \omega_\sigma) \geq f(\sigma, \omega) \quad \forall \omega \in \Omega \tag{10}$$

Let $f^*(\omega)$ be the optimal value of the flow time for the realization $\omega \in \Omega$. Then by definition, we also have

$$f(\sigma, \omega) \geq f^*(\omega) \quad \forall \omega \in \Omega \tag{11}$$

Using 10 and 11 we have,

$$f(\sigma, \omega_\sigma) \geq f^*(\omega) \quad \forall \omega \in \Omega \tag{12}$$

The above inequality implies that for any sequence $\sigma \in P$, the flow time corresponding to the worst case realization is greater than or equal to the optimal flow times for all realizations $\omega \in \Omega$. Since the above inequality holds for all $\sigma \in P$, it can be equivalently written as

$$\min_{\sigma \in P} f(\sigma, \omega_\sigma) \geq \max_{\omega \in \Omega} f^*(\omega) \tag{13}$$

Note that the left hand side of the above inequality is the objective of the ASMRSP. This proves the result.

In the past, significant work has has been done on developing approximation algorithms for $1|r_j|\sum C_j$ i.e. DSMSP with release times to minimize the total completion time of jobs. The best known approximation algorithm for $1|r_j|\sum C_j$ by Phillips et al. (1998) is a 2-approximation algorithm that produces non-preemeptive schedules from optimal preemptive schedules which can be easily determined using the *Shortest Remaining Processing Time* (SRPT) rule. It may be noted that for a given vector of release times and processing times, the optimal solution for $1|r_j|\sum C_j$ is also the optimal solution for $1|r_j|\sum C_j - r_j$. However the approximability of these two criteria may be very different as shown by Kellerer et al. (1999). Some of the reasonable approximation algorithms for $1|r_j|\sum C_j - r_j$ are the *Earliest Start Time* (EST) rule in which the shortest available job is assigned whenever the machine becomes free for assignment, or the *Earliest Completion Time* (ECT) rule in which the job with the earliest completion time (that may not be available yet) is assigned to the machine. Both the rules have a worst-case performance bound of $O(n)$. Kellerer et al. (1999) proposed an approximation algorithm with a sub-linear worst-case performance guarantee of $O(n^{1/2})$. They further showed that no constant ratio approximation algorithm can be expected for this problem by proving that there exists no polynomial time approximation algorithm with a worst-case performance bound of $O(n^{1/2-\epsilon})$, for any $\epsilon \geq 0$. It is clear that the bound obtained from the best known approximation algorithm is extremely weak for the problem under study in this paper. In the following section, we propose an exact method based on set-partitioning to solve the DSMSP.

**Exact Algorithm based on Set Partitioning**

Based on the preceding discussion, it is desirable to have an efficient exact algorithm to solve the deterministic single machine scheduling problem (DSMSP) to obtain a lower bound to the absolute robust single machine scheduling problem (ARSMSP). In this section, we propose an exact method based on set-partitioning to solve large instances of the DSMSP in small computation time. In this method, the set of all feasible assignments is generated apriori and is denoted by the set $L$. The assignment matrix is composed of the upper submatrix A and lower submatrix B. The upper submatrix A consists of $|L|$ columns and $|N|$ rows. In submatrix $A$, if column $l \in L$ represents the feasible assignment of job $i \in N$, then the entry in row $i$ is 1 while all other entries are zeroes. The lower submatrix B consists of $|L|$ columns and a single row for every discrete time interval in the planning horizon. Thus, in submatrix $B$, if column $l \in L$, represents the feasible assignment of job $i \in N$, then all entries corresponding to the time intervals in which the job $i$ is processed in the feasible assignment $l \in L$ are 1, while all the remaining entries are zeroes. To illustrate the procedure for the specific problem we are solving, consider the example containing two jobs, and four discrete time intervals in the planning horizon. Let us assume that both jobs have processing times of two time units, job 1 is released at time 1, while job 2 is released at the start of time 3, and hence can only be processed after that. Then the assignment matrix for the problem would look like as shown in Table 1. The first column represents the assignment of job 1 from time 1-2, and so on.

| | | | | |
|---|---|---|---|---|
| Job 1 | 1 | 1 | 1 | 0 |
| Job 2 | 0 | 0 | 0 | 1 |
| Time 1 | 1 | 0 | 0 | 0 |
| Time 2 | 1 | 1 | 0 | 0 |
| Time 3 | 0 | 1 | 1 | 1 |
| Time 4 | 0 | 0 | 1 | 1 |

Table 1: Assignment matrix for a simple example of set partitioning problem

We assume the following input data to be available for the set partitioning model:

$N =$ set of jobs
$H =$ set of discrete time intervals in the planning horizon
$L =$ set of feasible assignments
$t =$ $1, ..., |H|$ discrete time intervals in the planning horizon
$l =$ $1, ..., |L|$ feasible assignments
$d_l =$ delay associated with assignment $l$
$h_l =$ processing time associated with assignment $l$

The assignment matrix coefficients are defined as follows.

$$A_{il} = \begin{cases} 1 & \text{if the feasible assignment } l \text{ represents an assignment for job } i; \\ 0 & \text{otherwise.} \end{cases}$$

$$B_l^t = \begin{cases} 1 & \text{if job is being processed in time interval } t \text{ in assignment } l; \\ 0 & \text{otherwise.} \end{cases}$$

There is only a single decision variable for selection of feasible assignments in the optimal solution which is defined as follows.

$$\tau_l = \begin{cases} 1 & \text{if assignment } l \text{ is part of the optimal solution;} \\ 0 & \text{otherwise.} \end{cases}$$

The set partitioning model to solve the single machine scheduling problem with release times is formulated as shown below:

$$\min \sum_{j} (d_l \tau_l + h_l \tau_l) \tag{14}$$

$$\text{s.t.} \sum_{l} (A_{il} \tau_l) = 1 \quad \forall i \epsilon N \tag{15}$$

$$\sum_{l} (B_l^t \tau_l) \leq 1 \quad \forall t \epsilon H \tag{16}$$

$$\tau_l \epsilon \{0, 1\} \quad \forall l \epsilon L \tag{17}$$

In the above model, the objective (14) is to minimize the total flow time of the jobs, which includes the delays and the total processing times of the jobs. Note that the objective function can be equivalently expressed as the minimization of the sum of delays only, since the sum of processing times of the jobs given by $\sum_l (h_l \tau_l)$ is a constant. Thus in the proposed set partitioning model, the processing times are only used to build the matrix $B$. Constraints (15) ensure that each job must have exactly one feasible assignment in the optimal solution. Constraints (16) ensure that in a given time interval, at most one job can be processed. While the growth in the number of variables and constraints in the set-partitioning approach is much faster as compared to the mixed integer programming formulation discussed earlier, it can be used to solve the DSMSP to optimality almost instantaneously for even large problem size containing up to one hundred jobs, as validated by numerical experiments.

### 3.3.2 Robust Scheduling with release times

In the following, we discuss some properties of robust schedules with the performance criterion as the total flow time of the jobs. The release time of each job $i \in N$ lies in the interval $[\underline{r_i}, \overline{r_i}]$, and the processing time interval of job $i$ is $[\underline{p_i}, \overline{p_i}]$.

**ARSMSP with release times** The absolute robust single machine scheduling problem (ARSMSP) with release times can be mathematically formulated as (1)-(4) In Section 4, we describe two alternative heuristic algorithms to solve the ARSMSP with release times. The idea is to determine the sequence with the best worst-case performance by iteratively evaluating the worst-case performance of job sequences and moving towards better solutions using different search mechanisms. Clearly, in order to determine the sequence with the best worst-case absolute performance, it is essential to first formulate the problem of evaluating the worst case scenario for a given sequence $\sigma \in P$. Note that it is not straightforward to solve this problem by a simple enumeration technique, since the release times and processing times of all jobs are specified as independent ranges, thus implying an infinite number of possible realizations. Thus we have the following useful result that allows us to restrict our attention to only a subset of the realizations. To state the result, we define an *extreme point scenario* for a given sequence of jobs as one in which each job $i$ in the sequence assumes a release time $r_i$ equal to $\underline{r_i}$ or $\overline{r_i}$, and processing time $p_i$ equal to $\underline{p_i}$ or $\overline{p_i}$. The result is stated as follows

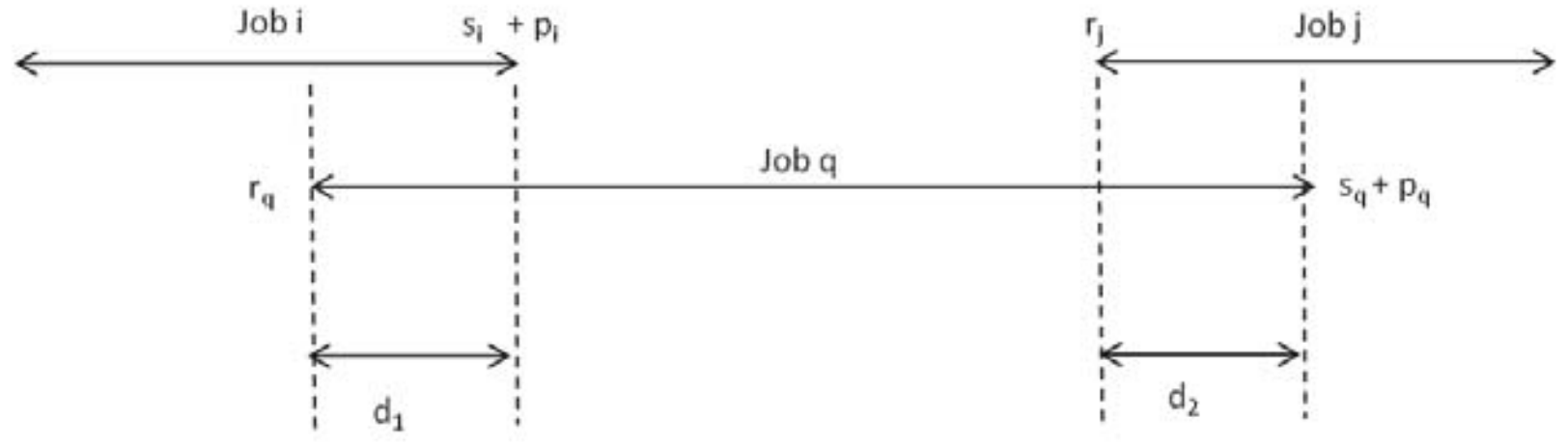

Figure 1: *Jobs* $i$, $q$ *and* $j$ *in an ordered sequence*

RESULT 4: *For the ARSMSP with* $n$ *jobs and uncertainty in both release times and processing times of jobs, there exists a worst-case scenario* $\omega_\sigma$ *for sequence* $\sigma \in P$*, that belongs to a subset of cardinality* $2^{n-2}$ *of the extreme point scenarios of* $\sigma$.

Proof: Consider an ordered sequence $N_\sigma$ of jobs, in which jobs $i$, $q$ and $j$ are consecutively ordered, that is, $i \prec q \prec j$. When job $q$ is the first or the last job in the sequence, jobs $i$ and $j$ respectively, may be considered as fake jobs. We assume that the release times and processing times of all jobs in the sequence except job $q$ are given (and unchangeable), and we want to show that there is an extreme point scenario of release time and processing time corresponding to job $q$, for which the job sequence assumes its worst case value. We define the following notations to illustrate the proof. Let $d_1$ be the overlap between the release time of job $q$ given by $r_q$ and the finishing time of processing job $i$ as given by $s_i + p_i$. Similarly let $d_2$ be the overlap between the end time of processing job $q$ given by $s_q + p_q$ and the release time of job $j$ given by $r_j$. This is graphically shown in Figure 1.

To obtain the worst case value, we need to maximize the sum of $d_1$ and $d_2$. We consider the following three cases:

- Case I: Job $q$ is the first job in the sequence. In this case, $d_1 = 0$ and $d_2 = \max(0, s_q + p_q - r_j)$. It is easy to see that there is a worst case scenario corresponding to $p_q = \overline{p_q}$ and $r_q = \overline{r_q}$.
- Case II: Job $q$ is the last job in the sequence. In this case, $d_1 = \max(0, s_i + p_i - r_q)$ and $d_2 = 0$. Again, it is easy to see that there is a worst case scenario corresponding to $r_q = \underline{r_q}$ and $p_q = \overline{p_q}$.
- Case III: When job $q$ lies somewhere in between, $d_1 = \max(0, s_i + p_i - r_q)$ and $d_2 = \max(0, s_q + p_q - r_j)$. By inspection, it can be inferred that $d_1 + d_2$ is maximized when $p_q = \overline{p_q}$ and $r_q = \overline{r_q}$ or $\underline{r_q}$.

Summarizing the above cases, there is a single unique extreme point scenario corresponding to the worst case contingency in cases I and II. For case III, for each of the $n$-2 possible positions of job $q$ in the sequence, there are 2 realizations of release times and a single realization of processing time for which the worst case value of the sequence can be obtained. Thus for $n$ jobs in a sequence, there exists a worst case scenario belonging to a subset of $2^{n-2}$ possible realizations. This proves the result.

The above result indicates that in order to determine the worst case scenario for a given sequence from the set of infinite possible realizations of release times and processing times, attention can be restricted to a subset of cardinality $2^{n-2}$ of extreme point scenarios. However this number can also be significantly large for large value of $n$. In the following, we show that the problem of finding the worst case realization for a given sequence can be formulated and solved as a mixed integer linear program (MILP).

In solving the absolute worst case performance problem (AWCPP) for a given sequence of jobs, we can preset the processing time of each job $i$ to its maximum value $\overline{p_i}$ using *Result 4*. Then the decision variables in the AWCPP include the start times of processing $s_i$ and the release times $r_i$ of jobs in the given sequence. The AWCPP for a given ordered sequence of jobs $N_\sigma$ can be stated as follows:

$$(AWCPP)\max\sum_{i\in N_\sigma}(s_i - r_i + \overline{p_i}) \tag{18}$$

$$s_1 = r_1 \tag{19}$$

$$s_i = \max\left(r_i, s_{i-1} + \overline{p_{i-1}}\right) \quad \forall i \in N_\sigma, i \geq 2 \tag{20}$$

$$r_i \in [\underline{r_i}, \overline{r_i}] \qquad \forall i \in N_\sigma \tag{21}$$

In the above model, the processing times constraints (19) state that the processing of the first job in the sequence starts as soon as it is released. The constraints (20) state that the processing of each subsequent job in the sequence should start as soon as the job is released and the processing of the previous job in the sequence has finished. The constraints (20) are not linear, but can be linearized using standard techniques (see Watters (1967)). To begin with we introduce two sets of additional variables $\eta_i$ and $\rho_i$ for all jobs $i \in N$. Then the constraints (20) can be equivalently expressed as

$$s_i = r_i + \eta_i \quad \forall i \in N_\sigma, i \geq 2 \tag{22}$$

$$s_i = s_{i-1} + \overline{p_{i-1}} + \rho_i \quad \forall i \in N_\sigma, i \geq 2 \tag{23}$$

$$\eta_i \rho_i = 0 \quad \forall i \in N_\sigma, i \geq 2 \tag{24}$$

To linearize constraints (24) we introduce binary variables $u_{ik}$ and $v_{ik}$ for all jobs $i \in N$, for a large enough positive integer $K$ such that $k \leq K$. Note that the product $\eta_i\rho_i$ is of the form $\sum_{t\leq K^2}\sum_{k\leq K}\sum_{k'\leq K} C_t u_{ik} v_{ik'}$, where the $C_t$ terms are constants. For the product $\eta_i\rho_i$ to be equal to zero, each term $C_t u_{ik} v_{ik'}$ should be equal to zero. This entails one or both the binary variables, $u_{ik}$ and $v_{ik'}$, to be equal to zero. This can be mathematically modeled as $u_{ik} + v_{ik'} \leq 1$. Thus we have the linearized version,

$$\eta_i = \sum_{k \le K} 2^k u_{ik} \qquad \forall i \in N_\sigma, i \ge 2 \tag{25}$$

$$\rho_i = \sum_{k \le K} 2^k v_{ik} \qquad \forall i \in N_\sigma, i \ge 2 \tag{26}$$

$$u_{ik} + v_{ik'} \le 1 \quad \forall i \in N_\sigma, i \ge 2, \forall k', k \le K \tag{27}$$

$$u_{ik}, v_{ik} \in \{0, 1\} \quad \forall i \in N_\sigma, i \ge 2, \forall k \le K \tag{28}$$

Following the above discussion, replacing $\eta_i$ and $\rho_i$ from constraints (25) - (26), the AWCPP can be rewritten as a mixed integer linear program as follows

$$(AWCPP) \max \sum_{i \in N_\sigma} (s_i - r_i + \overline{p_i}) \tag{29}$$

$$s_1 = r_1 \tag{30}$$

$$s_i = r_i + \sum_{k \le K} 2^k u_{ik} \qquad \forall i \in N_\sigma, i \ge 2 \tag{31}$$

$$s_i = s_{i-1} + \overline{p_{i-1}} + \sum_{k \le K} 2^k v_{ik} \qquad \forall i \in N_\sigma, i \ge 2 \tag{32}$$

$$u_{ik} + v_{ik'} \le 1 \quad \forall i \in N_\sigma, i \ge 2, \forall k', k \le K \tag{33}$$

$$u_{ik}, v_{ik} \in \{0, 1\} \quad \forall i \in N_\sigma, i \ge 2, \forall k \le K \tag{34}$$

$$r_i \in [\underline{r_i}, \overline{r_i}] \qquad \forall i \in N_\sigma \tag{35}$$

Thus given a sequence $\sigma \in P$, the worst case scenario can be determined by solving the above MILP. Note that in the above formulation, for $|N|$ jobs, the number of variables is of the order of $|N||K|$ and the number of constraints is of the order of $|N||K|^2$. From the computational experiments, the above MILP was found to be solvable almost instantaneously for even large problem size. The ARSMSP with release times given by 1-4 is solved using heuristic techniques described in the following section.

# 4 Solution Algorithms to the ARSMSP with Release Times

In this section, we present two alternative heuristic methods to obtain optimal or near-optimal solutions for the absolute robust single machine scheduling problem with uncertainty in release times and processing times.

## 4.1 Iterated Local Search

To begin with, we implement a simple heuristic based on iterated local search. In this method, we start with a random initial solution and perform a local search on the neighborhood of this sequence. In our implementation, the local search neighborhood $N_{LS}$ of a given sequence is defined as the set of sequences obtained by swapping two adjacent jobs in the original sequence. The performance measure for a given sequence is the worst-case scenario value of the sequence obtained by solving the AWCPP as given by 29-35. In case the local search improves the current solution, the local search solution is accepted as the new current solution and the local search is performed again. When the algorithm is stuck at a local minimum for too long, the algorithm is restarted with a new initial solution. The algorithm is terminated when the elapsed time from the beginning crosses a threshold computational time limit. The algorithm is described in Algorithm 1.

**Algorithm 1** Iterated Local Search Algorithm

**Require:** Set $N$ of jobs, set $M$ of scenarios
  Construct an initial feasible solution
  currentBestSolution ← initialSolution
  bestWorstCaseScenarioValue ← worstCaseScenarioValue(currentBestSolution)
  **while** timeLimit ≤ ilsTimeLimit **do**
    $x'$ = LocalSearch(currentBestSolution, $N_{LS}$)
    **if** worstCaseScenarioValue($x'$)< bestWorstCaseScenarioValue **then**
      bestWorstCaseScenarioValue ← worstCaseScenarioValue($x'$)
      currentBestSolution ← $x'$
    **end if**
    **if** solution value does not improve over time = timeRandomRestartILS **then**
      reinitialize currentSolution and start all over
    **end if**
  **end while**

## 4.2 Variable Neighborhood Search Algorithm

In this section, we propose the metaheuristic popularly known as the variable neighborhood search (VNS) in the literature. The algorithm was initially developed by Hansen and Mladenovic (1997). The main idea of the variable neighborhood search algorithm is to explore multiple neighborhood structures systematically instead of a single neighborhood, and escape local minima (in the case of

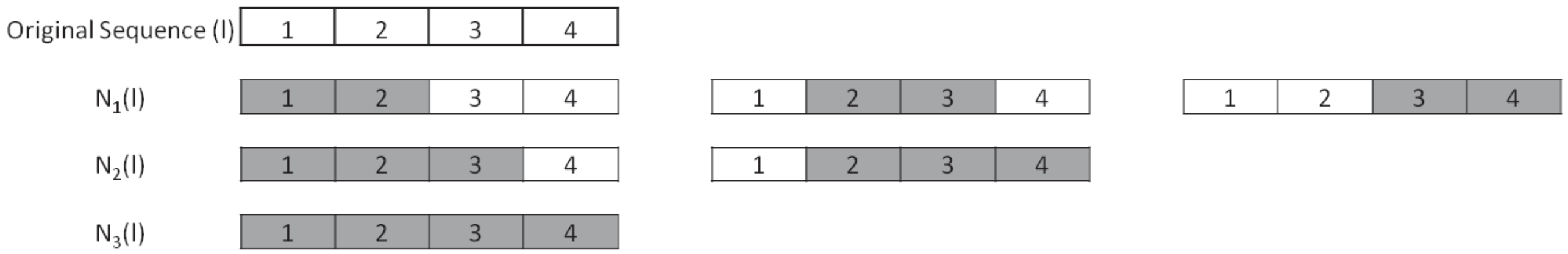

Figure 2: *VNS Neighborhood Structures for a given sequence 1-2-3-4*

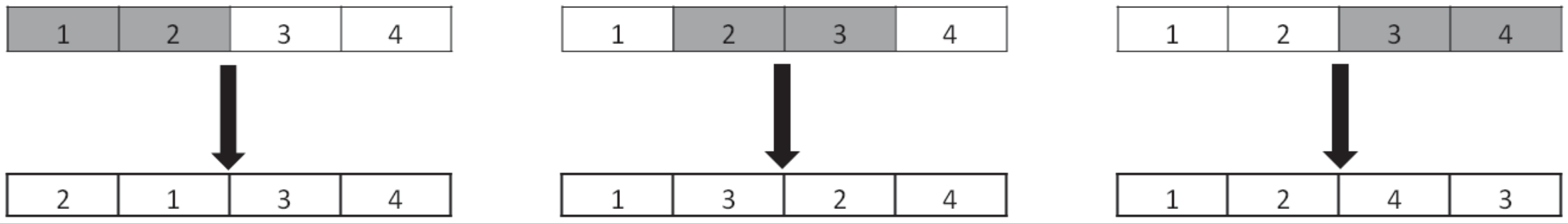

Figure 3: *VNS Neighborhood* $N_1$*(1-2-3-4)*

minimization). In our implementation of the method, the $k^{th}$ neighborhood structure, $N_k(\ell)$ of a given sequence $\ell$ is the set of sequences obtained by permuting the subset of jobs that are at most $k$ indices apart in the original sequence. It naturally follows that a sequence containing $n$ jobs has $n$-1 neighborhood structures. This is graphically represented in the Figure (2), where the permutable subset of jobs are shown in the blocks shaded in grey. Note that the neighborhoood structure $N_1(\ell)$ contains three candidate solutions as shown in Figure (3)

In the implementation of the VNS, we start with an initial feasible solution $x$. Iteratively starting from $k$=1, the *shaking procedure* is applied in which a random neighbor $x'$ is generated in the $N_k$ neighborhood of $x$. The shaking procedure is important as it prevents the algorithm from getting trapped at a local minimum. Thereafter a *local search* is carried out in the $N_{LS}$ neighborhood of $x'$, where the $N_{LS}$ neighborhood has a similar definition to the one described previously for the iterated local search method. The performance measure for a given sequence is the worst-case scenario value of the sequence obtained by solving the AWCPP as given by 29-35. If the local search solution $x''$ is found to be better than the current solution $x$, the search continues with the local search solution $x''$ as the new starting point, and $k$ is re-initialized to be equal to 1. If no improvement is found in the $N_k$ neighborhood, then $x$ remains the starting point for randomly generating a neighboring solution from the subsequent neighborhood $N_{k+1}$. When the current solution does not improve over a certain predefined time limit, the whole procedure is repeated starting from $k$=1 with a different initial solution. The algorithm is terminated when the time elapsed from the beginning crosses a threshold computational time limit. The algorithm is described in Algorithm 2:

**Algorithm 2** Variable Neighborhood Search Algorithm

```
Require: Set N of jobs, set M of scenarios
  Construct an initial feasible solution
  currentBestSolution ← initialSolution
  bestWorstCaseScenarioValue ← worstCaseScenarioValue(currentBestSolution)
  while timeLimit ≤ vnsTimeLimit do
    k=1
    while k ≤ (|N|-1) do
      Shaking Procedure
      x′ = GenerateNeighbor(currentBestSolution, N_k)
      Local Search
      x″ = LocalSearch(x′, N_LS)
      if worstCaseScenarioValue(x″)< bestWorstCaseScenarioValue then
        bestWorstCaseScenarioValue ← worstCaseScenarioValue(x″)
        currentBestSolution ← x″
        k=1
      else
        k++
      end if
      if solution value does not improve over time = timeRandomRestartVNS then
        reinitialize currentSolution and start all over
      end if
    end while
  end while
```

# 5 Computational Results and Analysis

## 5.1 Generation of Instances

The proposed heuristic algorithms were tested and validated through extensive numerical experiments based on artificial instances. The algorithms were implemented in JAVA programming language, and computational tests were run on an Intel Core i7 (2.80 GHz) processor and used a 32-bit version of CPLEX 12.2.

The experimental design adopted for the computational study consists of test problems involving $|N|$=7, 15, 20, 30 and 50 jobs and a single machine. For each problem size, 20 instances were tested. Based on the degree of stochastic variability in the release times and processing times of the jobs, the test instances are categorized into four different sets. For each category, the instances are generated by randomly drawing the lower and upper ends of the release time range and the processing time range of the jobs. The lower end of the release time range $\underline{r_i}$ is drawn from a uniform distribution of integers on the interval $\underline{r_i} \in [0, 5\mu]$ for four different values of $\mu$ ($\mu$=2,3,4 and 6). For $\mu = 2$ and 3, $\overline{r_i}$ is equal to $\underline{r_i} + 10$. On the other hand, for $\mu = 4$ and 6, $\overline{r_i}$ is equal to $\underline{r_i} + 20$. The lower end of the processing time range $\underline{p_i}$ is drawn from a uniform distribution of integers on the interval [1,4], while the upper end of the processing time range is equal to $\underline{p_i} + 6$. Five problem instances are tested for each combination of $|N|$ and $\mu$, resulting in a total of 100 problem instances.

## 5.2 Discussion of Results

The computational results obtained from the algorithms discussed previously are shown in the Tables (2)-(6). Using *Result 3* and *Result 4*, the lower bound is computed by taking the maximum of the optimal solution values over a subset of extreme point scenarios. The optimal solution for a given realization is obtained using the set partitioning algorithm described earlier.

For $|N|$=7 jobs, the optimal solution to the ARSMSP is calculated using an exhaustive search algorithm, and thus it is possible to determine the strength of the lower bound. As evident from Table (2), the lower bound is not too strong, and with increasing $\mu$ value, implying a larger uncertainty in the release times of the jobs, the bound weakens. For large problem size, it can be expected that the bound is even weaker.

From the results tables, it can be inferred that in general, the variable neighborhood search (VNS) algorithm is the more superior method to generate robust schedules. Based on a trial analysis, the computational time limit for test instances corresponding to a given combination of $|N|$ and $\mu$ was set to a certain value. It can be seen that for $|N|$=7, the VNS algorithm is able to generate optimal solutions for all instances in a computational time of few seconds. The iterated local search (ILS) method on the other hand is able to generate optimal solutions for close to 50% of the problem instances in the computational time limit of 100 seconds. For larger problem size with $|N|$ = 15, 20 30 and 50 jobs, the worst case value for a given sequence of jobs is determined by solving the mixed integer linear program 29-35 using K=10. On an average, the instances were found to converge faster for small $\mu$ value, that is, smaller uncertainty in the release times of the jobs. As can be seen from the results, the VNS and ILS algorithms converge to approximately the same solution for a few instances. Although the gap with respect to the lower bound is pretty large for most test instances, but since the bound is a weak one as established previously, it is difficult to comment on the absolute solution performance of the algorithms for these instances. Figure (4 ) shows the convergence of the solution for the test instance C18 over a computational time limit of 5 hours for the VNS and ILS methods. Note that the solution value may remain stable for a long time, before it begins to improve again.

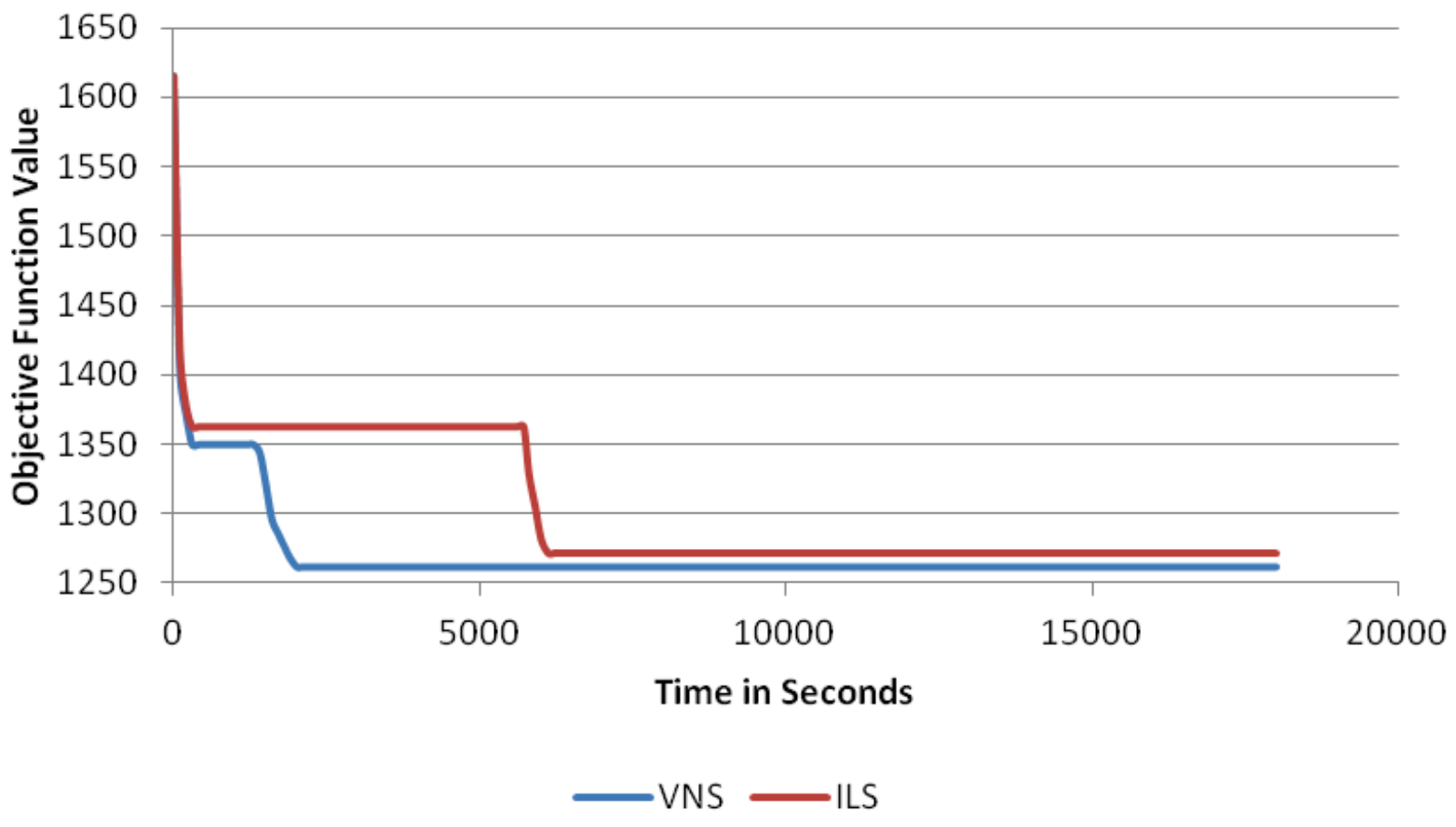

Figure 4: *Convergence of Instance C18 over a computational time limit of 5 hours*

From the computational experiments it was found that there is a certain degree of variance in the output solution values when a given problem instance was tested using a given algorithm. To study the behavior of the algorithms in more depth, we conduct a simulation study in which a test instance is run 50 times using a given algorithm and the resulting output solution values are plotted against the associated probability of finding a solution in the corresponding output range of values. The plots for some of the instances are shown in Figures (5)-(8). From the plots, the following observations can be made:

- In general, the mean of the output values for the VNS was found to be around the same or smaller than the ILS, implying that on an average, the VNS algorithm performs better than the ILS for most instances.

- There is a larger probability of finding a good solution using the VNS as compared to the ILS, as indicated by the frequency of the output solution values in the the low cost range as shown in the figures.

- The VNS is however less stable than the ILS as evident from the concentration of the output solution values in a single output range for the ILS, as represented by the peak in the distribution curve for the ILS.

Thus for a given instance, the VNS algorithm is expected to perform better on an average with a higher probability of finding a good solution, but there is also a larger variance in the output solution values returned by the VNS algorithm.

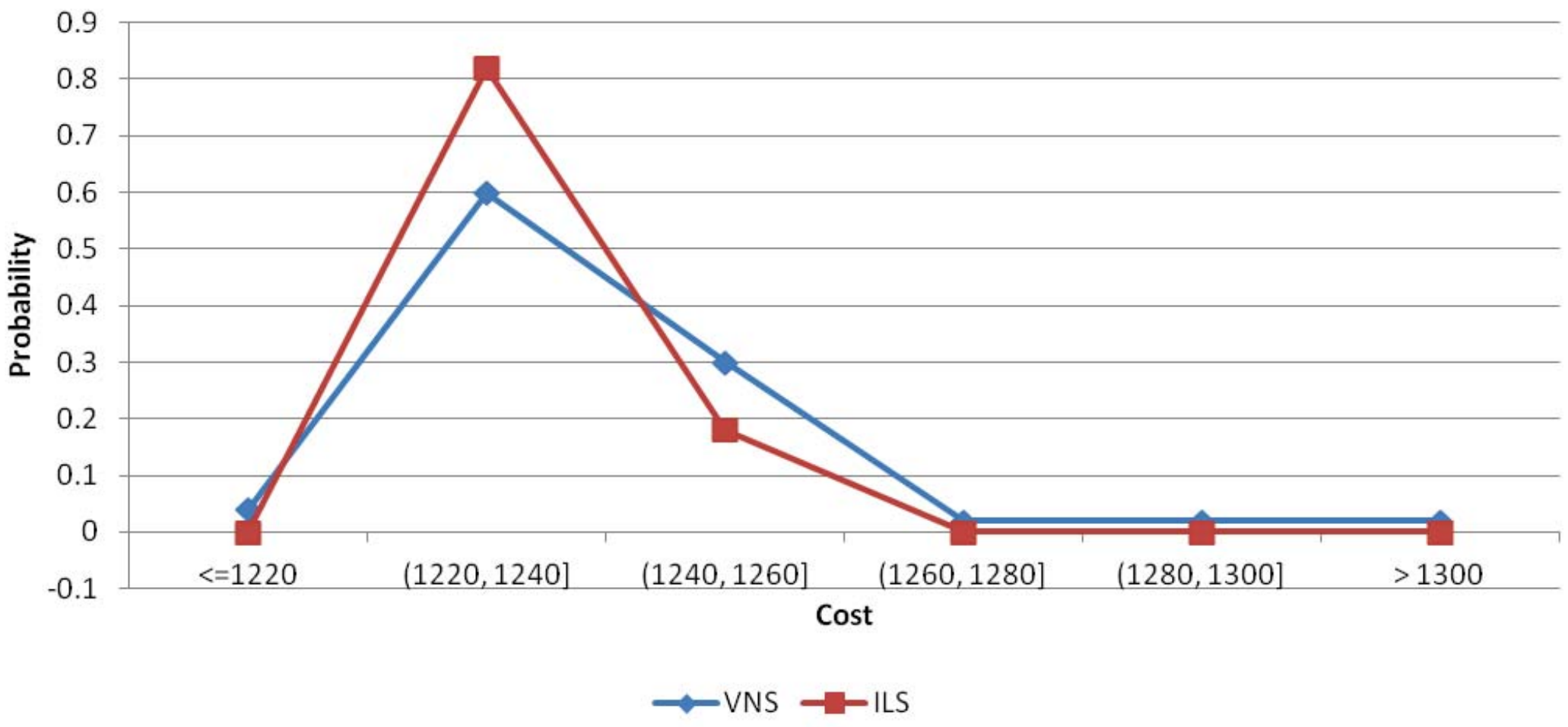

Figure 5: *Distribution of the output solution values for 50 simulation runs on instance C6 for a computational time limit of 300 seconds*

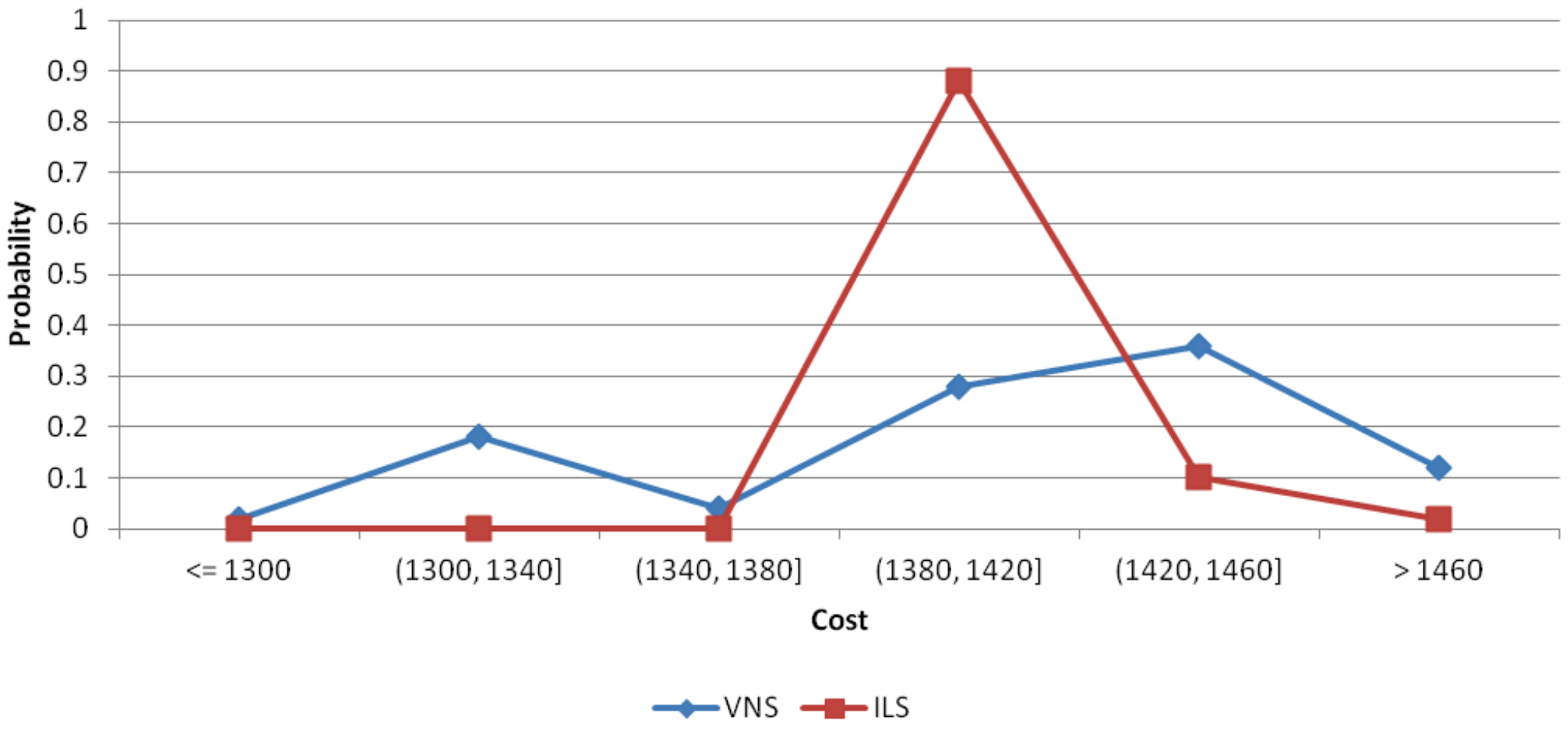

Figure 6: *Distribution of the output solution values for 50 simulation runs on instance C12 for a computational time limit of 300 seconds*

# 6 Conclusions and Future Work

This study demonstrates the complexity in dealing with uncertainty in release times and processing times of jobs in a proactive manner for the most basic form of the machine scheduling problem. In our problem, the release times and processing times of jobs are specified as independent ranges of values with unknown probability distributions. The performance criterion is the total flow time of all jobs and the robustness measure is the realized outcome for the worst-case contingency over the set of all possible scenarios.

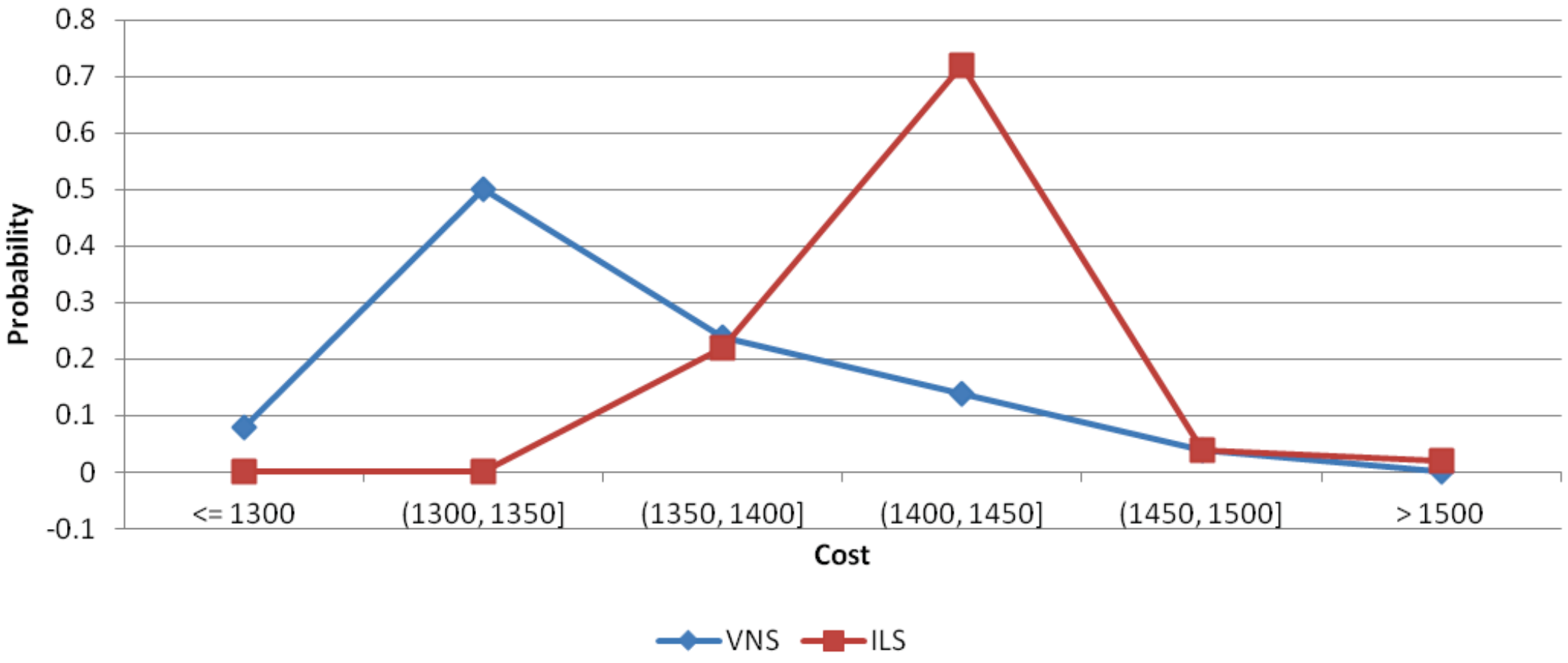

Figure 7: *Distribution of the output solution values for 50 simulation runs on instance C16 for a computational time limit of 300 seconds*

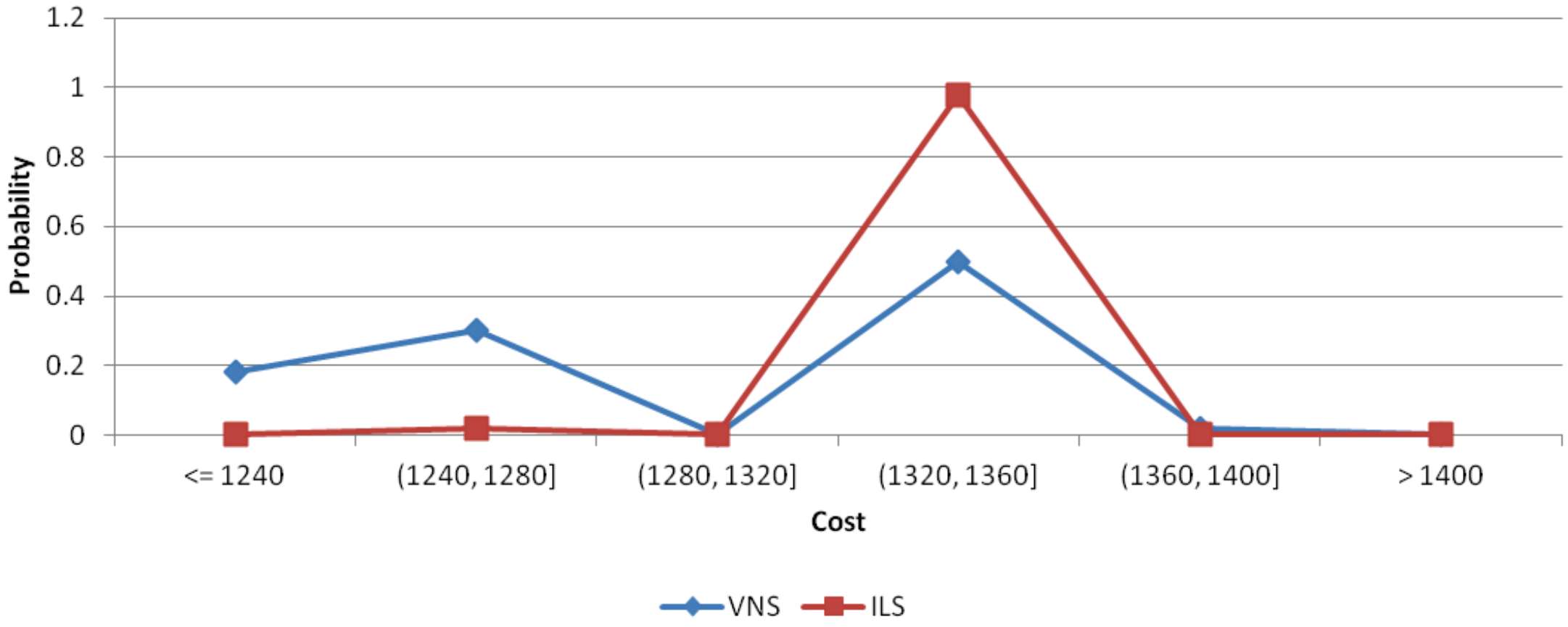

Figure 8: *Distribution of the output solution values for 50 simulation runs on instance C19 for a computational time limit of 300 seconds*

In previous research, the uncertainty in the release times of the jobs was largely ignored in the robust scheduling context. To the best of our knowledge, this is the first study that illustrates the added complexity in considering uncertainty in release times. We show that in order to solve the absolute robust single machine scheduling problem for $n$ jobs, we can restrict our attention to a subset of cardinality $2^n$ of the extreme point scenarios from the set of infinite possible realizations of release times and processing times. We propose heuristic algorithms based on variable neighborhood search (VNS) and iterated local search (ILS) to generate schedules with the best performance in the worst case contingency. The VNS algorithm was able to solve all instances with $|N|=7$ jobs to optimality. For larger problem size, on an average, the VNS was found to perform better than ILS with a larger associated probability of finding good solutions. However, the VNS was found to be less stable than the ILS as indicated by the variance in the output solution values.

As part of future work, the proposed methodology for the single machine scheduling problem can be extended to more than one machine. There is further scope for research on developing robust schedules for the single machine scheduling problem with different robustness measures such as the maximum regret or maximum relative deviation with respect to the corresponding optimal solution over the set of all possible realizations. There can also be other performance criteria such as the sum of completion times of all jobs or the total tardiness of all jobs beyond the specified due times for finishing.

Appendix A1

Table 2: Computational results for generated instances with |N|=7

| **Instance** | **Lower Bound**[1] | **Optimal Solution**[2] | | **VNS** | | **ILS** | | **% Gap**[3] |
|---|---|---|---|---|---|---|---|---|
| | | cost | time | cost | time[4] | cost | time[5] | |
| $\mu = 2$ | | | | | | | | |
| A1 | 188 | 212 | 4 | 212 | 5 | 212 | 2 | 11.32% |
| A2 | 189 | 204 | 4 | 204 | 7 | 217 | 100 | 7.35% |
| A3 | 192 | 200 | 4 | 200 | 13 | 207 | 100 | 4.00% |
| A4 | 168 | 180 | 5 | 180 | 7 | 192 | 100 | 6.67% |
| A5 | 189 | 198 | 4 | 198 | 4 | 203 | 100 | 4.55% |
| Mean | | | | | | | | **6.78%** |
| $\mu = 3$ | | | | | | | | |
| A6 | 174 | 189 | 5 | 189 | 2 | 189 | 1 | 7.94% |
| A7 | 126 | 137 | 4 | 137 | 4 | 139 | 100 | 8.03% |
| A8 | 160 | 183 | 5 | 183 | 7 | 183 | 2 | 12.57% |
| A9 | 162 | 202 | 4 | 202 | 3 | 202 | 1 | 19.80% |
| A10 | 194 | 200 | 4 | 200 | 4 | 206 | 100 | 3.00% |
| Mean | | | | | | | | **10.27%** |
| $\mu = 4$ | | | | | | | | |
| A11 | 135 | 173 | 5 | 173 | 7 | 175 | 100 | 21.97% |
| A12 | 146 | 154 | 4 | 154 | 10 | 154 | 76 | 5.19% |
| A13 | 120 | 133 | 4 | 133 | 5 | 133 | 4 | 9.77% |
| A14 | 180 | 209 | 4 | 209 | 83 | 213 | 100 | 13.88% |
| A15 | 121 | 151 | 5 | 151 | 13 | 151 | 8 | 19.87% |
| Mean | | | | | | | | **14.14%** |
| $\mu = 6$ | | | | | | | | |
| A16 | 148 | 178 | 4 | 178 | 8 | 178 | 5 | 16.85% |
| A17 | 141 | 189 | 4 | 189 | 9 | 189 | 5 | 25.40% |
| A18 | 150 | 192 | 4 | 192 | 8 | 200 | 100 | 21.88% |
| A19 | 134 | 171 | 4 | 171 | 12 | 173 | 100 | 21.64% |
| A20 | 156 | 183 | 5 | 183 | 7 | 194 | 100 | 14.75% |
| Mean | | | | | | | | **20.10%** |

[1]Using *Result 3* and *Result 4*, the lower bound is the maximum optimal value over a subset of extreme point scenarios.
[2]The optimal solution is determined using an exhaustive search algorithm.
[3]The gap indicates the optimality gap of the lower bound with respect to the optimal solution.
[4]A computational time limit of 100 seconds was set for all instances with |N|=7 jobs.
[5]A computational time limit of 100 seconds was set for all instances with |N|=7 jobs.

Appendix A2

Table 3: Computational results for generated instances with |N|=15

| Instance | Lower Bound | VNS | | ILS | |
|---|---|---|---|---|---|
| | | cost | time[6] | cost | time[7] |
| $\mu = 2$ | | | | | |
| B1 | 726 | 761 | 300 | 761 | 300 |
| B2 | 703 | 783 | 300 | 783 | 300 |
| B3 | 748 | 783 | 300 | 783 | 300 |
| B4 | 674 | 710 | 300 | 710 | 300 |
| B5 | 713 | 750 | 300 | 809 | 300 |
| $\mu = 3$ | | | | | |
| B6 | 716 | 757 | 300 | 757 | 300 |
| B7 | 654 | 690 | 300 | 747 | 300 |
| B8 | 712 | 742 | 300 | 742 | 300 |
| B9 | 630 | 657 | 300 | 665 | 300 |
| B10 | 593 | 624 | 300 | 637 | 300 |
| $\mu = 4$ | | | | | |
| B11 | 599 | 701 | 600 | 723 | 600 |
| B12 | 587 | 681 | 600 | 763 | 600 |
| B13 | 601 | 698 | 600 | 713 | 600 |
| B14 | 687 | 701 | 600 | 770 | 600 |
| B15 | 620 | 676 | 600 | 778 | 600 |
| $\mu = 6$ | | | | | |
| B16 | 671 | 763 | 600 | 763 | 600 |
| B17 | 663 | 723 | 600 | 723 | 600 |
| B18 | 675 | 748 | 600 | 772 | 600 |
| B19 | 689 | 771 | 600 | 841 | 600 |
| B20 | 724 | 871 | 600 | 889 | 600 |

[6]The computational time limit determined from a trial based analysis.
[7]The computational time limit determined from a trial based analysis.

Appendix A3

Table 4: Computational results for generated instances with |N|=20

| Instance | Lower Bound | VNS | | ILS | |
|---|---|---|---|---|---|
| | | cost | time[8] | cost | time[9] |
| $\mu = 2$ | | | | | |
| C1 | 1254 | 1322 | 600 | 1322 | 600 |
| C2 | 1305 | 1410 | 600 | 1410 | 600 |
| C3 | 1289 | 1369 | 600 | 1434 | 600 |
| C4 | 1259 | 1395 | 600 | 1399 | 600 |
| C5 | 1259 | 1338 | 600 | 1338 | 600 |
| $\mu = 3$ | | | | | |
| C6 | 1117 | 1228 | 600 | 1228 | 600 |
| C7 | 1226 | 1274 | 600 | 1364 | 600 |
| C8 | 1237 | 1317 | 600 | 1365 | 600 |
| C9 | 1206 | 1328 | 600 | 1328 | 600 |
| C10 | 1144 | 1269 | 600 | 1356 | 600 |
| $\mu = 4$ | | | | | |
| C11 | 1131 | 1342 | 900 | 1329 | 900 |
| C12 | 1094 | 1369 | 900 | 1390 | 900 |
| C13 | 1208 | 1362 | 900 | 1363 | 900 |
| C14 | 1084 | 1312 | 900 | 1312 | 900 |
| C15 | 1130 | 1351 | 900 | 1426 | 900 |
| $\mu = 6$ | | | | | |
| C16 | 1063 | 1213 | 900 | 1225 | 900 |
| C17 | 1038 | 1144 | 900 | 1259 | 900 |
| C18 | 1085 | 1355 | 900 | 1363 | 900 |
| C19 | 1067 | 1215 | 900 | 1337 | 900 |
| C20 | 1116 | 1357 | 900 | 1357 | 900 |

[8] The computational time limit determined from a trial based analysis.

[9] The computational time limit determined from a trial based analysis.

Appendix A4

Table 5: Computational results for generated instances with |N|=30

| Instance | Lower Bound | VNS | | ILS | |
|---|---|---|---|---|---|
| | | cost | time[10] | cost | time[11] |
| $\mu = 2$ | | | | | |
| D1 | 2865 | 3252 | 600 | 3286 | 600 |
| D2 | 2595 | 2985 | 600 | 3126 | 600 |
| D3 | 2768 | 3063 | 600 | 3120 | 600 |
| D4 | 2647 | 3095 | 600 | 2982 | 600 |
| D5 | 2799 | 3328 | 600 | 3335 | 600 |
| $\mu = 3$ | | | | | |
| D6 | 2526 | 2871 | 600 | 2786 | 600 |
| D7 | 2452 | 3170 | 600 | 3183 | 600 |
| D8 | 2525 | 2946 | 600 | 3160 | 600 |
| D9 | 2673 | 3147 | 600 | 3154 | 600 |
| D10 | 2809 | 3263 | 600 | 3326 | 600 |
| $\mu = 4$ | | | | | |
| D11 | 2378 | 2870 | 900 | 3085 | 900 |
| D12 | 2543 | 3224 | 900 | 3320 | 900 |
| D13 | 2294 | 2612 | 900 | 2593 | 900 |
| D14 | 2551 | 3135 | 900 | 3193 | 900 |
| D15 | 2362 | 2826 | 900 | 2835 | 900 |
| $\mu = 6$ | | | | | |
| D16 | 2296 | 3004 | 900 | 3092 | 900 |
| D17 | 2224 | 2836 | 900 | 2748 | 900 |
| D18 | 2205 | 3262 | 900 | 3320 | 900 |
| D19 | 2190 | 3082 | 900 | 3256 | 900 |
| D20 | 2153 | 2892 | 900 | 2914 | 900 |

[10] The computational time limit determined from a trial based analysis.
[11] The computational time limit determined from a trial based analysis.

Appendix A5

Table 6: Computational results for generated instances with |N|=50

| Instance | Lower Bound | VNS | | ILS | |
|---|---|---|---|---|---|
| | | cost | time[12] | cost | time[13] |
| $\mu = 2$ | | | | | |
| E1 | 7761 | 9326 | 900 | 9363 | 900 |
| E2 | 7608 | 9060 | 900 | 9245 | 900 |
| E3 | 7650 | 9001 | 900 | 8989 | 900 |
| E4 | 7655 | 8670 | 900 | 8910 | 900 |
| E5 | 7650 | 9568 | 900 | 9574 | 900 |
| $\mu = 3$ | | | | | |
| E6 | 7423 | 8718 | 900 | 9127 | 900 |
| E7 | 7057 | 8129 | 900 | 8675 | 900 |
| E8 | 7337 | 8616 | 900 | 8642 | 900 |
| E9 | 7846 | 9111 | 900 | 9185 | 900 |
| E10 | 7135 | 8710 | 900 | 8722 | 900 |
| $\mu = 4$ | | | | | |
| E11 | 7647 | 9500 | 1200 | 9767 | 1200 |
| E12 | 7672 | 9642 | 1200 | 9655 | 1200 |
| E13 | 7099 | 8836 | 1200 | 8848 | 1200 |
| E14 | 7478 | 8972 | 1200 | 9207 | 1200 |
| E15 | 7354 | 9126 | 1200 | 9321 | 1200 |
| $\mu = 6$ | | | | | |
| E16 | 6861 | 8598 | 1200 | 8596 | 1200 |
| E17 | 7062 | 9105 | 1200 | 9116 | 1200 |
| E18 | 7309 | 8662 | 1200 | 8671 | 1200 |
| E19 | 6984 | 8677 | 1200 | 8671 | 1200 |
| E20 | 7076 | 8630 | 1200 | 8957 | 1200 |

[12]The computational time limit determined from a trial based analysis.
[13]The computational time limit determined from a trial based analysis.